\documentclass [12pt,leqno]{article}

\newcounter{conjecture}
\newcounter{remark}
\newcounter{corollary}

\newcommand{\eqnsection}{
   \renewcommand{\theequation}{\thesection.\arabic{equation}}
   \makeatletter
   \csname @addtoreset\endcsname{equation}{section}
   \makeatother}
\newtheorem{theorem}{Theorem}
\newtheorem{lemma}{Lemma}

\def \be{\begin{equation}}
\def \ee{\end{equation}}
\def \bt{\begin{theorem}}
\def \et{\end{theorem}}
\def \bea{\begin{eqnarray}}
\def \eea{\end{eqnarray}}
\def \bas{\begin{eqnarray*}}
\def \eas{\end{eqnarray*}}

\def \al{\alpha}
\def \bb{\beta}
\def \ga{\gamma}

\def \de{\delta}
\def \De{\Delta}
\def \ep{\epsilon}

\def \la{\lambda}

\def \ff{\infty}
\def \wh{\widehat}
\def \wt{\widetilde}

\def \rar{\rightarrow}

\def \({\left(}
\def \){\right)}

\def \lc{\left\{}
\def \rc{\right\}}

\def \nn{\nonumber}

\def \bc{\begin{center} }
\def \ec{\end{center} }
\def \bs{\begin{slide} }
\def \es{\end{slide} }

\def\square{{\vcenter{\vbox{\hrule height.3pt
        \hbox{\vrule width.3pt height5pt \kern5pt
           \vrule width.3pt}
        \hrule height.3pt}}}}
\def\qed{{\hfill $\square$ \bigskip}}

\eqnsection
\begin{document}

\def\wh{\widehat}
\def\ol{\overline}

\title{A stochastic calculus proof of the  CLT  for the $L^{2}$ modulus of continuity of  local time }

\author{  Jay Rosen\thanks
     {This research was  supported, in part, by grants from the National Science
Foundation and PSC-CUNY.}}


\maketitle

\bibliographystyle{amsplain}

\begin{abstract} We give a stochastic calculus proof of the Central Limit Theorem
\[ { \int ( L^{ x+h}_{t}- L^{ x}_{ t})^{ 2}\,dx- 4ht\over h^{ 3/2}}
\stackrel{\mathcal{L}}{\Longrightarrow}c\(\int ( L^{ x}_{ t})^{ 2}\,dx\)^{1/2}\,\,\eta\]
as $h\rar 0$ for Brownian local time $L^{ x}_{ t}$. Here $\eta$  is an independent   normal random variable with mean zero and variance
one.
\end{abstract}



\section{Introduction}

 In   \cite{lp}  we obtain  almost sure limits for the
$L^{ p}$ moduli  of continuity of  local times of a very wide class of
symmetric L\'evy processes. More specifically, if $\{L^{ x }_{ t}\,;\,(x,t)\in R^{ 1}\times  R^{  1}_{ +}\}$ denotes Brownian local time then  for all
$  p\ge 1$, and all
$t\in R_+$
\be
\lim_{ h\downarrow 0}  \int_{a}^{ b} \bigg|{  L^{ x+h}_{ t} -L^{ x }_{
t}\over\sqrt{h}}\bigg|^p\,dx =2^pE(|\eta|^p)
\int_a^b |L^{ x }_{ t}|^{ p/2}\,dx\label{as.1}
\ee for all
$a,b
$ in the extended real line  almost surely, and also  in $L^m$, $m\ge 1$.
(Here $\eta$  is    normal random variable with mean zero and variance
one.) In particular when $p=2$ we have
\begin{equation}
\lim_{ h\downarrow 0}  \int  {  (L^{ x+h}_{ t} -L^{ x }_{
t})^{ 2}\over h}\,dx =4t, \hspace{.2 in}\mbox{ almost surely.} \label{rp3.1}
\end{equation}
We refer to $\int    (L^{ x+h}_{ t} -L^{ x }_{
t})^{ 2} \,dx$ as the  $L^{2}$ modulus of continuity of Brownian  local time.

In our recent paper \cite{CLMR} we obtain  the central limit theorem corresponding to (\ref{rp3.1}). 

\begin{theorem}\label{theo-clt2} For each fixed $t$
\begin{equation} { \int ( L^{ x+h}_{t}- L^{ x}_{ t})^{ 2}\,dx- 4ht\over h^{ 3/2}}
\stackrel{\mathcal{L}}{\Longrightarrow}c\(\int ( L^{ x}_{ t})^{ 2}\,dx\)^{1/2}\,\,\eta\label{5.0weak}
\end{equation}  as $h\rar 0$, with $c=\({64 \over 3}\)^{ 1/2}$. Equivalently
\begin{equation} { \int ( L^{ x+1}_{t}- L^{ x}_{ t})^{ 2}\,dx- 4t\over t^{ 3/4}}
\stackrel{\mathcal{L}}{\Longrightarrow}c\(\int ( L^{ x}_{ 1})^{ 2}\,dx\)^{1/2}\,\,\eta\label{5.0tweak}
\end{equation}
as $t\rar\ff$. Here $\eta$  is an independent   normal random variable with mean zero and variance
one.
\end{theorem} 

It can be shown that
\begin{equation}
E\(\int ( L^{ x+1}_{ t}- L^{ x}_{ t})^{ 2}\,dx\)=4\( t-{2t^{ 1/2} \over \sqrt{2\pi }  }\)+O( 1).\label{9.13}
\end{equation}
so that (\ref{5.0tweak}) can be written as
\begin{equation} { \int ( L^{ x+1}_{t}- L^{ x}_{ t})^{ 2}\,dx- E\(\int ( L^{ x+1}_{ t}- L^{ x}_{ t})^{ 2}\,dx\)\over t^{ 3/4}}
\stackrel{\mathcal{L}}{\Longrightarrow}c\(\int ( L^{ x}_{ 1})^{ 2}\,dx\)^{1/2}\,\,\eta\label{5.0tweake}
\end{equation}
with a similar statement for (\ref{5.0weak}).

Our proof of Theorem \ref{theo-clt2} in \cite{CLMR} is rather long and involved. We use the method of moments, but rather than study the asymptotics of the moments of (\ref{5.0weak}), which seem intractable, we study the moments of the  analogous expression where the fixed time $t$ is replaced by an independent exponential time of mean $1/\la$. An important part of the proof is then to `invert the Laplace transform' to obtain the asymptotics of the moments for fixed $t$.

 The purpose of this paper is to give a new and shorter proof of Theorem \ref{theo-clt2} using stochastic integrals, following the approach of \cite{Yor, YW}. Our proof makes use of certain differentiability properties of the double and triple intersection local time, $ \al_{2,t}(x)$ and $\al_{3,t}(x,y)$, which are formally given by 
   \begin{equation}
 \al_{2,t}(x)=\int_{0}^{t} \int_{0}^{s} \de (W_{s}-W_{r}-x ) \,dr \,ds \label{6.3}
 \end{equation}
 and
\begin{equation}
\al_{3,t}(x,y)= \int_{0}^{t} \int_{0}^{s}\int_{0}^{r} \de (W_{r}-W_{r'}-x) \de (W_{s}-W_{r}-y ) \,dr'  \,dr \,ds.\label{6.14}
\end{equation}
More precisely, let $f(x)$ be a smooth positive symmetric function with compact support and 
 $\int f(x)\,dx=1$. Set $f_{\ep}(x)={1 \over \ep}f(x/\ep)$. Then 
   \begin{equation}
 \al_{2,t}(x)=\lim_{\ep\rar 0}\int_{0}^{t} \int_{0}^{s} f_{\ep} (W_{s}-W_{r}-x ) \,dr \,ds \label{6.3}
 \end{equation}
  and
\bea
&&
\al_{3,t}(x,y)\nonumber\\
&&= \lim_{\ep\rar 0}\int_{0}^{t} \int_{0}^{s}\int_{0}^{r} f_{\ep} (W_{r}-W_{r'}-x) f_{\ep} (W_{s}-W_{r}-y ) \,dr'  \,dr \,ds\label{6.14}
\eea
 exist almost surely and in all $L^{p}$, are independent of the particular choice of $f$, and are continuous in $(x,y,t)$ almost surely, \cite{jc}.  It is easy to show, see \cite[Theorem 2]{djcrilt}, that for any measurable $\phi(x)$
\bea
&&
  \int_{0}^{t} \int_{0}^{s} \phi (W_{s}-W_{r}) \,dr \,ds= \int \phi(x)\al_{2,t}(x)\,dx\label{6.14r2}
\eea
and for any measurable $\phi(x,y)$
\bea
&&
  \int_{0}^{t} \int_{0}^{s}\int_{0}^{r} \phi (W_{r}-W_{r'},\,W_{s}-W_{r}) \,dr'  \,dr \,ds\nonumber\\
&&= \int \phi(x,y)\al_{3,t}(x,y)\,dx\,dy.\label{6.14r}
\eea
 
To express the differentiability properties of $\al_{2,t}(x)$ and $\al_{3,t}(x,y)$ which we need, let us set
\begin{equation}
v(x)=\int_{0}^{\ff}e^{-s/2}p_{s}(x)\,ds=e^{-|x|}.\label{6.14s}
\end{equation}
The following result is   \cite[Thorem 1]{djcrilt}.
\begin{theorem}\label{lem-diff}
 \begin{equation}
 \ga_{2,t}(x)=:\al_{2,t}(x)-tv(x)\label{6.18b}
 \end{equation}
 and
 \begin{equation}
 \ga_{3,t}(x,y)=:\al_{3,t}(x,y)-\ga_{2,t}(x)v(y)-\ga_{2,t}(y)v(x)-tv(x)v(y)\label{6.18}
 \end{equation}
 are $C^{1}$ in $(x,y)$ and $\nabla \ga_{2,t}(x), \nabla \ga_{3,t}(x,y)$ are continuous in $(x,y,t)$.
 \end{theorem} 
 
Our new proof of Theorem \ref{theo-clt2} is given in Section  \ref{sec-stoch}.

 Our original motivation for studying the
asymptotics of
$\int   (L^{ x+h}_{ t} -L^{ x }_{
t})^{ 2}\,dx$ comes from our interest in the Hamiltonian
\begin{equation} H_{ n}=\sum_{ i,j=1,\,i\neq j}^{ n}1_{ \{S_{ i}=S_{  j} \}}-{1
\over 2}\sum_{ i,j=1,\,i\neq j}^{ n}1_{
\{|S_{ i}-S_{ j}|=1 \}},\label{rp5c.4}
\end{equation}
  for the  critical attractive 
random polymer in dimension one,
\cite{HK}, where $\{S_{ n} \,;\,n=0,1,2,\ldots\}$ is a simple random walk on
$Z^{ 1}$. Note that
$ H_{ n}=\sum_{x\in Z^{ 1}}\(l_{ n}^{ x}-l_{ n}^{ x+1}\)^{ 2}$, where $l_{ n}^{
x}=\sum_{ i=1}^{ n}1_{ \{S_{ i}=x \}}$   is the local time for the  random walk
$S_{ n}$.

 \section{A stochastic calculus approach}\label{sec-stoch}

 By \cite[Lemma 2.4.1]{book} we have that
    \begin{equation}
L^{ x}_{ t}=\lim_{\ep\rar 0}\int_{0}^{t}  f_{\ep} (W_{s} -x )  \,ds \label{6.3l}
 \end{equation}
 almost surely, with convergence locally uniform in $x$. Hence
  \begin{eqnarray}
 &&  \int   L^{ x+h}_{ t}L^{ x}_{ t} \,dx\nonumber\\
 &&=  \int \lim_{\ep\rar 0} \( \int_{0}^{t} f_{\ep} (W_{s}-(x+h))\,ds  \)
 \( \int_{0}^{t} f_{\ep} (W_{r}-x)\,dr  \)\,dx
 \label{6.2}\\
  &&=  \lim_{\ep\rar 0}\int  \( \int_{0}^{t} f_{\ep} (W_{s}-(x+h))\,ds  \)
 \( \int_{0}^{t} f_{\ep} (W_{r}-x)\,dr  \)\,dx
 \nonumber\\
 &&=  \lim_{\ep\rar 0}\int_{0}^{t} \int_{0}^{t} f_{\ep}\ast f_{\ep} (W_{s}-W_{r}-h ) \,dr \,ds    \nonumber\\
 &&=  \lim_{\ep\rar 0}\int_{0}^{t} \int_{0}^{s} f_{\ep}\ast f_{\ep} (W_{s}-W_{r}-h ) \,dr \,ds\nonumber\\
 &&\hspace{1.5 in}+ \lim_{\ep\rar 0} \int_{0}^{t} \int_{0}^{r}f_{\ep}\ast f_{\ep} (W_{r}-W_{s}+h )\,ds  \,dr    \nonumber\\
 &&=\al_{2,t}(h)+\al_{2,t}(-h).\nn
 \end{eqnarray} 
 Note that
 \begin{equation}
 \int ( L^{ x+h}_{ t}- L^{ x}_{ t})^{ 2}\,dx=2\( \int   (L^{ x}_{ t} )^{ 2}\,dx-
  \int   L^{ x+h}_{ t}L^{ x}_{ t} \,dx\)\label{6.1}
 \end{equation}
 and thus
 \begin{equation}
 \int ( L^{ x+h}_{ t}- L^{ x}_{ t})^{ 2}\,dx=2\( 2\al_{2,t}(0)- \al_{2,t}(h)- \al_{2,t}(-h)\).\label{6.4}
 \end{equation}
 Hence we can prove Theorem \ref{theo-clt2} by showing that  for each fixed $t$
\begin{equation} { 2\( 2\al_{2,t}(0)- \al_{2,t}(h)- \al_{2,t}(-h)\)- 4ht\over h^{ 3/2}}
\stackrel{\mathcal{L}}{\Longrightarrow}c\sqrt{\al_{2,t}(0)}\,\,\eta\label{5.0weaksi}
\end{equation}  
as $h\rar 0$, with $c=\({128 \over 3}\)^{ 1/2}$. Here we used the fact, which follows from (\ref{6.2}), that $\int ( L^{ x}_{ 1})^{ 2}\,dx= 2\,\,\al_{2,t}(0)$.

In proving (\ref{5.0weaksi}) we will need the following Lemma. Compare Tanaka's formula,  \cite[Chapter VI, Theorem 1.2]{RY}.
\begin{lemma}\label{lem-Ito}For any $a\in R^{1}$
\begin{eqnarray}
&&
 \al_{2,t}(a)=2\int_{0}^{t}  (W_{t}-W_{s}-a)^{+}\,ds -2\int_{0}^{t}  (W_{0}-W_{s}-a)^{+}\,ds \nn\\
&& \hspace{1.5 in}-2(-a)^{+}t-2 \int_{0}^{t} \int_{0}^{s} 1_{\{W_{s}-W_{r}>a \}}\,dr\,dW_{s}. \label{it.1}
\end{eqnarray}
\end{lemma}

{\bf  Proof of Lemma \ref{lem-Ito}: }Set
\begin{equation}
g_{\ep}(x)=\int_{0}^{\ff} yf_{\ep}(x-y)\,dy\label{it.2}
\end{equation}
so that
\begin{equation}
g'_{\ep}(x)=\int_{0}^{\ff} yf'_{\ep}(x-y)\,dy=\int_{0}^{\ff} f_{\ep}(x-y)\,dy\label{it.3}
\end{equation}
and consequently
\begin{equation}
g''_{\ep}(x)=f_{\ep}(x). \label{it.4}
\end{equation}
 
 Let
 \begin{equation}
F_{a}(t,x)=  \int_{0}^{t} g_{\ep}(x-W_{s}-a)\,ds. 
\label{6.5}
 \end{equation} 
 Then by Ito's formula applied to $F_{a}(t,W_{t}) $ we have 
 \begin{eqnarray}
 &&\int_{0}^{t} g_{\ep}(W_{t}-W_{s}-a)\,ds -\int_{0}^{t} g_{\ep}(W_{0}-W_{s}-a)\,ds 
 \label{6.6}\\
 && = \int_{0}^{t} g_{\ep}(-a)\,ds+ \int_{0}^{t} \int_{0}^{s} g'_{\ep}(W_{s}-W_{r} -a)\,dr\,dW_{s}\nonumber\\
 &&\hspace{1 in}+{1 \over 2}\int_{0}^{t} \int_{0}^{s} g''_{\ep}(W_{s}-W_{r} -a)\,dr\,ds.\nn
 \end{eqnarray}
 It is easy to check that locally uniformly
 \begin{equation}
\lim_{\ep\rar 0}g_{\ep}(x)=x^{+}\label{it.5}
\end{equation}
and hence using (\ref{it.4}) we obtain 
\begin{eqnarray}
&&
 \al_{2,t}(a)=2\int_{0}^{t}  (W_{t}-W_{s}-a)^{+}\,ds -2\int_{0}^{t}  (W_{0}-W_{s}-a)^{+}\,ds \nn\\
&& \hspace{1.2 in}-2(-a)^{+}t-2\lim_{\ep\rar 0}\int_{0}^{t} \int_{0}^{s} g'_{\ep}(W_{s}-W_{r} -a)\,dr\,dW_{s}.  \label{it.1a}
\end{eqnarray}
From (\ref{it.3}) we can see that $\sup_{x} |g'_{\ep}(x)|\leq 1$ and 
 \begin{equation}
\lim_{\ep\rar 0}g'_{\ep}(x)=1_{\{x>0\}}+{1 \over 2}1_{\{x=0\}}.\label{it.5a}
\end{equation}
Thus by the dominated convergence theorem
\begin{equation}
\lim_{\ep\rar 0}\int_{0}^{t}E\(\(  \int_{0}^{s} \lc g'_{\ep}(W_{s}-W_{r} -a)-  1_{\{W_{s}-W_{r}>a \}}   \rc\,dr     \)^{2}\)\,ds =0\label{it.6}
\end{equation}
which completes the proof of our Lemma.\qed

If we now set 
\begin{eqnarray}
 &&
J_{h}(x)=2x^{+}-(x-h)^{+}-(x+h)^{+}\nn\\ 
&&\hspace{.45in}  =\left\{\begin{array}{ll}
-x-h&\mbox{ if }-h\leq x\leq 0\\
x-h&\mbox{ if }0\leq x\leq h.
\end{array}
\right.\label{6.8a}
 \end{eqnarray}
 and 
 \begin{eqnarray}
 K_{h}(x) &=&21_{\{x>0\}}-1_{\{x>h\}}-1_{\{x>-h\}}\label{6.9}\\
 &=&1_{\{0<x\leq h\}}-1_{\{-h<x\leq 0\}}\nn
 \end{eqnarray}
we  see from Lemma \ref{lem-Ito} that 
\begin{eqnarray}
&&2\lc 2\al_{t}(0)- \al_{t}(h)- \al_{t}(-h)\rc-4ht
\label{6.10}\\
&&= 4 \int_{0}^{t} J_{h}(W_{t}-W_{s})\,ds -4\int_{0}^{t} J_{h}(W_{0}-W_{s})\,ds  \nonumber\\
&&\hspace{1 in}-4\int_{0}^{t} \int_{0}^{s} K_{h}(W_{s}-W_{r} )\,dr\,dW_{s}.\nn
\end{eqnarray}
By (\ref{6.8a})
 \begin{equation}
 \int_{0}^{t} J_{h}(W_{t}-W_{s})\,ds= \int  J_{h}(W_{t}-x)L^{x}_{t}\,dx=O(h^{2}\sup_{x}L^{x}_{t})\label{6.8b}
 \end{equation}
 and similarly for $\int_{0}^{t} J_{h}(W_{0}-W_{s})\,ds$. Hence to prove (\ref{5.0weaksi}) it suffices to show that  for each fixed $t$
\begin{equation} { \int_{0}^{t} \int_{0}^{s} K_{h}(W_{s}-W_{r} )\,dr\,dW_{s}\over h^{ 3/2}}
\stackrel{\mathcal{L}}{\Longrightarrow}\({8 \over 3}\)^{ 1/2}\sqrt{\al_{2,t}(0)}\,\,\eta\label{5.0weaksia}
\end{equation}  
as $h\rar 0$. Let
\begin{equation}
M^{h}_{t}=h^{- 3/2}\int_{0}^{t} \int_{0}^{s} K_{h}(W_{s}-W_{r} )\,dr\,dW_{s}.\label{6.11}
\end{equation}
It follows from the proof of Theorem 2.6 in \cite[Chapter XIII]{RY},
(the Theorem of Papanicolaou, Stroock, and Varadhan) 
that to establish (\ref{5.0weaksia})
it suffices to show that
  \begin{equation}
\lim_{h\rar 0} \langle M^{h}, W\rangle_{t} =0\label{6.26h1}
 \end{equation}
 and 
  \begin{equation}
\lim_{h\rar 0} \langle M^{h}, M^{h}\rangle_{t} ={8 \over 3} \al_{2,t}(0)\label{6.26h}
 \end{equation}
 uniformly in $t$ on compact intervals.

By (\ref{6.14r2}), and using the fact that $ K_{h}(x) = K_{1}(x/h) $, we have that
\begin{eqnarray}
 \langle M^{h}, W\rangle_{t}&=&h^{- 3/2}\int_{0}^{t} \int_{0}^{s} K_{h}(W_{s}-W_{r} )\,dr\,ds
\label{6.14r5}\\
&=&h^{- 3/2}\int K_{h}(x)\al_{2,t}(x)\,dx   \nonumber\\
&=&h^{- 1/2}\int K_{1}(x)\al_{2,t}(hx)\,dx   \nonumber\\
&=&\int_{0}^{1}{ \al_{2,t}(hx)- \al_{2,t}(-hx) \over h^{1/2}}\,dx.   \nonumber
\end{eqnarray}
But $v(hx)- v(-hx)$, so by Lemma  \ref{lem-diff} we have that
\begin{eqnarray}
&&\al_{2,t}(hx)- \al_{2,t}(-hx)=\ga_{2,t}(hx)- \ga_{2,t}(-hx)=O(h) \label{6.14r6}
\end{eqnarray} 
  which completes the proof of (\ref{6.26h1}).
  
We next analyze 
\begin{eqnarray} 
&&
 \langle M^{h}, M^{h}\rangle_{t}=h^{-3}\int_{0}^{t}\( \int_{0}^{s} K_{h}(W_{s}-W_{r} )\,dr\)^{2}\,ds \label{6.12}\\
&&=h^{-3}\int_{0}^{t}\( \int_{0}^{s} K_{h}(W_{s}-W_{r} )\,dr\)\( \int_{0}^{s} K_{h}(W_{s}-W_{r'} )\,dr'\) \,ds \nn\\
&&=h^{-3}\int_{0}^{t}\( \int_{0}^{s}\int_{0}^{r}  K_{h}(W_{s}-W_{r'} )K_{h}(W_{s}-W_{r} )\,dr'  \,dr\) \,ds \nn\\
&&+h^{-3}\int_{0}^{t}\( \int_{0}^{s}\int_{0}^{r'}  K_{h}(W_{s}-W_{r} )K_{h}(W_{s}-W_{r'} )\,dr  \,dr'\) \,ds. \nn
\end{eqnarray}
By (\ref{6.14r}) we have that
\begin{eqnarray}
&&\int_{0}^{t} \int_{0}^{s}\int_{0}^{r}  K_{h}(W_{s}-W_{r'} )K_{h}(W_{s}-W_{r} )\,dr'  \,dr \,ds
\label{6.13}\\
&&=  \int_{0}^{t} \int_{0}^{s}\int_{0}^{r}  K_{h}(W_{s}-W_{r} +W_{r}-W_{r'})K_{h}(W_{s}-W_{r} )\,dr'  \,dr \,ds \nonumber\\
&&=  \int  \int   K_{h}(x+y)K_{h}(y)\al_{3,t}(x,y)  \,dx \,dy. \nonumber
\end{eqnarray}
 Using $ K_{h}(x) = K_{1}(x/h) $ we have 
\begin{eqnarray}
&&h^{-3}\int_{0}^{t} \int_{0}^{s}\int_{0}^{r}  K_{h}(W_{s}-W_{r'} )K_{h}(W_{s}-W_{r} )\,dr'  \,dr \,ds
\label{6.15}\\
&&= h^{-3} \int  \int   K_{h}(x+y)K_{h}(y)\al_{3,t}(x,y)  \,dx \,dy \nonumber \\
&&= h^{-1} \int  \int  K_{1}(x+y) K_{1}(y)\al_{3,t}(hx,hy)  \,dx \,dy \nonumber\\
&&=  h^{-1}\int  \int   K_{1}(x)K_{1}( y)\al_{3,t}(h(x-y),hy)  \,dx \,dy \nonumber\\
&&= h^{-1} \int_{0}^{1}  \int_{0}^{1}  A_{3,t}(h,x,y)  \,dx \,dy \nonumber
\end{eqnarray}
where 
\begin{eqnarray}
&&
A_{3,t}(h,x,y)=\al_{3,t}(h(x-y),hy)-\al_{3,t}(h(-x-y),hy)\label{6.16}\\
&&\hspace{1 in}-\al_{3,t}(h(x+y),-hy)+\al_{3,t}(-h(x-y),-hy).\nn
\end{eqnarray}
 It remains to consider
 \begin{equation}
 \lim_{h\rar 0}{ A_{3,t}(h,x,y)\over h}.\label{6.17}
 \end{equation}

 We now use Lemma \ref{lem-diff}.
Using the fact that  $\ga_{3,t}(x,y), \ga_{2,t}(x)$  are continuously differentiable
 \begin{eqnarray}
 && \ga_{3,t}(h(x-y),hy)- \ga_{3,t}(h(-x-y),hy) 
 \label{6.19}\\
 && =h(x-y){\partial  \over \partial x}\ga_{3,t}(0,hy)-h(-x-y){\partial  \over \partial x} \ga_{3,t}(0,hy) +o(h)\nonumber\\
 && =2hx{\partial  \over \partial x} \ga_{3,t}(0,0) +o(h)\nonumber
 \end{eqnarray}
 and similarly
  \begin{eqnarray}
 && \ga_{3,t}(-h(x-y),-hy)- \ga_{3,t}(h(x+y),-hy) 
 \label{6.19}\\
 && =-h(x-y){\partial  \over \partial x} \ga_{3,t}(0,-hy)-h(x+y){\partial  \over \partial x} \ga_{3,t}(0,-hy)  +o(h)\nonumber\\
 && =-2hx{\partial  \over \partial x} \ga_{3,t}(0,0) +o(h)\nonumber
 \end{eqnarray}
 and these two terms cancel up to $o(h)$.
 
 Next,
 \begin{eqnarray}
 && \ga_{2,t}(h(x-y))v(hy)- \ga_{2,t}(h(-x-y))v(hy)
 \label{6.20}\\
 &&  \hspace{1 in}+ \ga_{2,t}(-h(x-y))v(-hy)-\ga_{2,t}(h(x+y))v(-hy) \nonumber\\
 && = h(x-y)\ga'_{2,t}(0)v(0)-h(-x-y) \ga'_{2,t}(0)v(0)\nonumber\\
 &&  \hspace{.6 in}-h(x-y) \ga'_{2,t}(0)v(0)-h(x+y)\ga'_{2,t}(0)v(0)+o(h) \nonumber\\
 &&  =o(h). \nonumber
 \end{eqnarray}
 
 On the other hand, using $v(x)=e^{-|x|}=1-|x|+  O(x^{2})$ we have 
  \begin{eqnarray}
 && v(h(x-y))\ga_{2,t}(hy)- v(h(-x-y))\ga_{2,t}(hy)
 \label{6.21}\\
 &&  \hspace{1 in}+ v(-h(x-y))\ga_{2,t}(-hy)-v(h(x+y))\ga_{2,t}(-hy) \nonumber\\
 && =-| h(x-y)|\ga_{2,t}(0)+|h(-x-y)|\ga_{2,t}(0)\nonumber\\
 &&  \hspace{.6 in}-|h(x-y)|\ga_{2,t}(0)+|h(x+y)|\ga_{2,t}(0)+o(h) \nonumber\\
 &&  =2h(|x+y|-|x-y|)\ga_{2,t}(0)+o(h). \nonumber
 \end{eqnarray}
 and similarly
   \begin{eqnarray}
  && v(h(x-y))v(hy)- v(h(-x-y))v(hy)
 \label{6.22}\\
 &&  \hspace{1 in}+ v(-h(x-y))v(-hy)-v(h(x+y))v(-hy) \nonumber\\
 && =-| h(x-y)|v(0)+|h(-x-y)|v(0)\nonumber\\
 &&  \hspace{.6 in}-|h(x-y)|v(0)+|h(x+y)|v(0)+O(h^{2}) \nonumber\\
 &&  =2h(|x+y|-|x-y|)v(0)+O(h^{2}). \nonumber
 \end{eqnarray}
 Putting this all together and using the fact that $\al_{2,t}(0)=\ga_{2,t}(0)+tv(0)$ we see that
   \begin{eqnarray}
  &&\int_{0}^{1}  \int_{0}^{1}  A_{3,t}(h,x,y)  \,dx \,dy \label{6.23}\\
  &&
  =2h\al_{2,t}(0)\int_{0}^{1}  \int_{0}^{1}(|x+y|-|x-y|))\,dx \,dy+o(h).\nn
 \end{eqnarray}
 Of course \begin{eqnarray}
 &&\int_{0}^{1}  \int_{0}^{1}(|x+y|-|x-y|))\,dx \,dy
 \label{6.24}\\
 && =\int_{0}^{1}  \int_{0}^{x}2y\,dy \,dx+\int_{0}^{1}  \int_{0}^{y}2x\,dx \,dy={2 \over 3}  \nonumber
 \end{eqnarray}
 so that 
 \begin{equation}
\lim_{h\rar 0}{\int_{0}^{1}  \int_{0}^{1}  A_{3,t}(h,x,y)  \,dx \,dy  \over h}={4 \over 3} \al_{2,t}(0).\label{6.25}
 \end{equation}
By (\ref{6.12}) this gives (\ref{6.26h}).\qed

 \def\noopsort#1{} \def\printfirst#1#2{#1}
\def\singleletter#1{#1}
      \def\switchargs#1#2{#2#1}
\def\bibsameauth{\leavevmode\vrule height .1ex
      depth 0pt width 2.3em\relax\,}
\makeatletter
\renewcommand{\@biblabel}[1]{\hfill#1.}\makeatother
\newcommand{\bysame}{\leavevmode\hbox to3em{\hrulefill}\,}

\bigskip
\noindent
\begin{tabular}{lll} 
      & Jay Rosen \\
      & Department of Mathematics \\
     &College of Staten Island, CUNY \\
     &Staten Island, NY 10314 \\ &jrosen30@optimum.net  
\end{tabular}


\begin{thebibliography}{10}


\bibitem{CLMR} X.~Chen, W. Li, M. Marcus and  J.~Rosen,  {\em   A CLT  for the $L^{2}$ modulus of continuity of local times of Brownian motion},\,   {\it Ann. Probab.},\, to appear. \,http://arxiv.org/pdf/0901.1102v1.pdf
 

   \bibitem{HK} R. van der Hofstad, A. Klenke and W. Konig, {\em The critical attractive 
random polymer in dimension one}. Jour. Stat. Phys. 106:3-4, 477-520 (2002)


 \bibitem{book} M. Marcus and J.~Rosen, {\em Markov Processes, Gaussian
Processes and Local Times}, Cambridge studies in advanced mathematics,
100, Cambridge University Press, Cambridge, England, 2006.

 \bibitem{lp} M. B. Marcus and J.~Rosen,
        {\em $L^{ p}$ moduli of continuity of Gaussian processes and
 local times
 of symmetric L\'evy processes},  Annals of Probab., to appear.
 
  \bibitem{RY} D. Revuz and M.~Yor,, {\em Continuous martingales and Brownian motion}, Springer, Berlin, 1998.

\bibitem{jc} J.~Rosen,  {\em Joint continuity of renormalized intersection local
times.}
 Ann. Inst. Henri Poincare,\,{\bf 32}\, (1996), 671--700. 
 
 \bibitem{dsilt} J.~Rosen, {\em  Derivatives of self-intersection local times.}
S\'{e}minaire de Probabilit\'{e}s,\,\,XXXVIII,\,  Springer-Verlag, New York ,
(2005), LNM 1857, 171-184. 


 \bibitem{djcrilt} J.~Rosen, {\em  Continuous Differentiability of Renormalized Intersection Local Times in $R^{1}$}, preprint.
 
 \bibitem{Yor} M.~Yor, {\em  Le drap brownien comme limite en lois des temps locaux lin\'{e}aires.}
S\'{e}minaire de Probabilit\'{e}s,\,\, XVII,\,  Springer-Verlag, New York ,
(1983), LNM 986, 89-105.

 \bibitem{YW} S.~Weinryb and M.~Yor, {\em  Le mouvement brownien de L\'{e}vy index\'{e} par $R^{3}$ comme limite centrale des temps locaux d'intersection.}
S\'{e}minaire de Probabilit\'{e}s,\,\, XXII,\,  Springer-Verlag, New York ,
(1988), LNM 1321, 225-248.


\end{thebibliography}
\end{document}